\documentclass[a4paper,12pt]{article}
\usepackage{etex}
\usepackage[utf8x]{inputenc}
\usepackage{caption, subcaption}
\usepackage{hyperref}
\usepackage{relsize}
\usepackage[all]{xy}
\usepackage{amsmath,amssymb,amsfonts,amscd, graphicx, latexsym, verbatim, multirow, color}
\usepackage{epsfig}
\usepackage{geometry} % to change the page dimensions
\usepackage{adjustbox}

\newcommand{\mtz}{{\mathbf M}_2({\Z})}

\newcommand{\mtmz}{{\mathbf M}_2^{n.s.}({\Z})}

\newcommand{\mtboxz}{{\mathbf M}_2^\Box({\Z})}
\newcommand{\mtr}{{\mathbf M}_2({\R})}

\newcommand{\Q}{\mathbf Q}
\newcommand{\Z}{\mathbf Z}
\newcommand{\pdet}{\mathrm{pdet}}

\newcommand{\inv}{^{-1}}
 \newcommand{\R}{\mathbf R}

\newcommand{\psl}{\mathrm{PSL}_2(\Z)}
\newcommand{\pgl}{\mathrm{PGL}_2(\Z)}
\newcommand{\pglz}{\mathrm{PGL}_2(\Z)}

\newcommand{\glq}{\mathrm{GL}_2(\Q)}
\newcommand{\pglq}{\mathrm{PGL}_2(\Q)}
\newcommand{\pslq}{\mathrm{PSL}_2(\Q)}
\newcommand{\PGL}{\mathrm{PGL}}

\newtheorem{corollary}{Corollary}

\newtheorem{proposition}{Proposition}
\usepackage{color, ulem}

\usepackage{tikz}
\usetikzlibrary{matrix,arrows}

\begin{document}

%opening
\title{A presentation of  PGL(2,Q)}
\author{A. Muhammed Uluda\u{g}\footnote{
{Department of Mathematics, Galatasaray University}
{\c{C}{\i}ra\u{g}an Cad. No. 36, 34349 Be\c{s}ikta\c{s}}
{\.{I}stanbul, Turkey}}}

\maketitle

\begin{abstract}
We give a conjectural presentation of the infinitely generated group PGL(2,Q)
with an infinite list of relators.
\end{abstract}

\section{Introduction}
The group $\pglq$ is defined as the projectivization of the group $\glq$. 
Our aim is to give a conjectural presentation of it. 

This group contains the subgroup $\pslq$ of $2\times 2$ rational projective matrices of determinant 1. 
One has the  inclusions
$$
\psl<\pgl<\pslq^\pm <\pglq < {\mathrm{PGL}_2(\R)},
$$
$$
\psl<\pslq <\pglq^+< {\mathrm{PSL}_2(\R)}, \mbox{ where }
$$

\noindent
$\pslq^\pm$ is the group of $2\times 2$ rational projective matrices of determinant $\pm$1. \\
$\pgl$ is  the group of $2\times 2$ integral projective matrices of determinant 1. \\
$\psl$ is  the group of $2\times 2$ integral projective matrices of determinant 1. \\
$\pglq^+$ is the group of $2\times 2$ rational matrices of determinant $>0$.\\

Every inclusion in the above list is of infinite index, except the first one. 
The modular group $\psl$ and its $\Z/2\Z$-extension $\pgl$ in this hierarchy are finitely generated and their presentations are known. Other groups 
are not finitely  generated and to our knowledge their presentations are not known. 

Dresden \cite{dresden} classifies the finite subgroups and finite-order elements (see Section \ref{finite} below) of this group.
Besides a few papers about some $\pglq$-cocycles related to some partial zeta values and to generalized  Dedekind sums,
(see \cite{solomon}, \cite{sczech}) we were unable to spot any works about this group.

An infinitely presented group is usually not a very friendly object, nevertheless, due to its connection to the lattice $\Z^2$, the group $\pglq$ merits a better treatment.  It turns out that its presentation is not so complicated, and one has the following parallelism between $\pglz$ and 
$\pglq$. The Borel subgroup $B(\Z)$ of $\pglz$, which by definition is the set of upper triangular elements, is generated by the translation $Tz:=1+z$ and the reflection $Vz:=-z$ (where we take the liberty to consider the elements of $\pglq$ as elements of the M\"obius group of $\mathbb P^1(\R)$). Since $T=KV$, the Borel subgroup is also generated by the involutions $V$ and $Kz:=1-z$, showing that $B(\Z)$ is the infinite dihedral group. The group $\pglz$ itself is generated by its  Borel subgroup $B(\Z)$ and the involution $Uz:=1/z$. Note that the derived subgroup of $B(\Z)$ is 
$\Z$. 

Likewise, $\pglq$ is generated by its Borel subgroup $B(\Q)$ and $U$. 
Here $B(\Q)$ is infinitely generated but nevertheless is quite similar to the infinite dihedral group in that its derived subgroup is $\Q$.
It admits the presentation
$$
B(\Q)\simeq \langle K, H_{p}\,|\, H_p^{-1}T^pH_p=T, \, [H_p, H_q]=1, \, p,q=-1,2,3,5,7\dots \rangle 
$$
where the elements $H_nz:=nz$ are the homotheties  (note that $T=KV$ and $V=H_{-1}$). Its abelianization is $\Q^\times$. These claims are easily verified by using the matrix description of $B(\Q)$.

The subgroup $\langle T, H_p \rangle<B(\Q)$ admits the presentation 
$$
\langle T, H_p \,|\, H_p^{-1}T^pH_p=T \rangle,
$$
i.e. it is the Baumslag-Solitar group $BS(p,1)$.
Its abelianization is $\Z\times\Z/(p-1)\Z $. See \cite{levitt}, \cite{epstein} for more about $BS(p,1)$.

The derived subgroup of $\pglq$ is $\pslq$. It is simple \cite{alperin}. 
The abelianization of $\pglq$ is the multiplicative group of integers  modulo square integers.
This latter is an infinitely generated 2-torsion group. 
It follows that any quotient of $\pglq$ is an abelian $2$-torsion group. 

\section{The monoid of integral matrices}
Let $\Lambda$ be the free abelian group of rank 2 and consider the monoid of all $\Z$-module morphisms  
$\Lambda \longrightarrow \Lambda$. 
By using the standard generators for $\Lambda$, we may identify this monoid with the 
monoid of $2\times 2$ matrices over $\Z$, denoted by $\mtz$. Thus we consider the monoid
$$
\mtz:=\left\{ \left(\begin{matrix} p&q \\ r&s \end{matrix}\right) \, \left|\right. p,q,r,s\in \Z\right\},
$$
under the usual matrix product.
Set-wise, $\mtz$ is just $\Z^4$. It admits the set of scalar matrices $\langle nI \,|\,n\in \Z\rangle$ as a submonoid. The quotient of $\mtz$ by this submonoid will be referred to as its {\it projectivization} and denoted by  $\mathbb P\mtz$.
Denote by $\mtmz$ the set of non-singular matrices and by $\mtboxz$ those with a square determinant.

The set of invertible elements inside $\mathbb P\mtz$  consists of projective two-by-two non-singular matrices with rational (or equivalently integral) entries and is denoted by $\pglq$.
To be more precise, denote by $[M]$ the projectivization of $M$.
Then 
$$
\pglq:=\bigl\{[M] \,: \, M\in \mtz, \quad \det(M)\neq 0\bigr\} =\mathbb P\mtmz.
$$
The group $\pglq$ is isomorphic to the group of integral  M\"obius transformations of $\mathbb P^1(\R)$.
It contains the projectivization of the set of matrices with rational entries and of determinant 1 as a subgroup, the group $\pslq$.
Hence,
$$
\pslq:=\bigl\{[M] \,: \, M\in \mtz, \quad \det(M) \mbox{ is a square}\bigr\}=\mathbb P\mtboxz.
$$
The multiplicative group $\Z^\times$ admits a subgroup which consists of square rationals
which we denote as $\Z^{2\times}$.
The group $\Z^\times/\Z^{2\times}$ of square-free rationals is the torsion abelian group generated by $\{x_{-1}\}\cup \{x_p: p \mbox{ prime}\}$
subject to the relations $x_i^2=1$ for each $i=-1, 2, 3, \dots$. It is isomorphic to $\Q^\times/\Q^{2\times}$.

The map $\det:\mtz\to \Z$ is equivariant under the action of $a\in \Z$ by
$$
\left(\begin{matrix} p&q \\ r&s \end{matrix}\right)\in \mtz\to 
\left(\begin{matrix} ap&aq \\ ar&as \end{matrix}\right)\in \mtz
$$
on the left and by $t\to a^{2}t$ on the right. Hence we may projectivize it as follows:

\begin{displaymath}
    \xymatrix{
    1\ar@{->}[r]&\langle\pm I\rangle\ar@{->}[r] \ar@{->}[d]&I\Z^\times \ar@{->>}[r]\ar@{->}[d]&\Z^{2\times}\ar@{->}[d]\\ 
    1\ar@{->}[r]&\mtboxz\ar@{->}[r] \ar@{->>}[d]& \mtmz\ar@{->>}[r]^\det\ar@{->>}[d]&\Z^\times\ar@{->>}[d]\\
    1\ar@{->}[r]&\pslq\ar@{->}[r] &\pglq  \ar@{->>}[r]^{{\pdet}} &\Z^\times/\Z^{2\times}\\
    }
\end{displaymath}

%Note that the derived subgroup of $\pglq$ is $\pslq$. Indeed $\pslq$ is included in the derived subgroup since it is the kernel of $\pdet$, which is a map onto the abelian group $\Q^\times/\Q^{2\times}$. 
%The derived subgroup is included in $\pslq$ since $\pdet([M,N])=1$ for every $M, N$.

\section{Generators}
\paragraph{The modular group.} It is well known (see \cite{coxeter}) that the extended modular group $\pgl$ admits the presentation
$$
\langle V,U,K\, |\, V^2=U^2=K^2=(VU)^2=(KU)^3=1 \rangle.
$$
The modular group is the subgroup $\psl$ is generated by $S(z)=-1/z$ and $L(z)= (z-1)/z$ and admits the presentation
$
\psl =\langle S, L \,|\, S^2=L^3=1\rangle.
$
\paragraph{Homotheties.}For $r=m/n\in \Q$, consider the homothety $H_r:z\mapsto rz$. In matrix form
$$
H_r=  \left[\begin{matrix} r&0 \\ 0&1 \end{matrix}\right] = \left[\begin{matrix} m&0 \\ 0&n \end{matrix}\right]
\implies H_r\inv=H_{1/r}, \quad H_rH_s=H_sH_r=H_{rs}.
$$
Hence every homothety is a product of ``prime" homotheties.
$H_{-1}$ is involutive. The remaining homotheties forms a  free abelian group of infinite rank. 
\paragraph{Involutions.}Let $r=m/n\in \Q$ and consider the involutions $I_r(z):=r/z=m/nz\in \pglq$, which is represented by the traceless matrix
$$
I_r=\left[\begin{matrix} 0&r \\ 1&0 \end{matrix}\right]=\left[\begin{matrix} 0&m \\ n&0 \end{matrix}\right].
$$
Note that $\pdet(I_p)=x_p$.  One has
$$
 I_rI_s=H_{r/s}, \quad I_sI_r=H_{s/r}, \quad [I_r, I_s]=I_rI_sI_r\inv I_s\inv= H_{r^2/s^2}\in \pslq
$$
In particular, $I_r$ and $I_s$ do not commute unless $r=\pm s$.

The set of involutions $\{I_r\,:\, r\in \Q\} $  generate a group containing the homotheties and which  consists of  two lines inside $\mathbb P\mtr$: 
$$
\langle I_r\, |\, r\in |\Q^\times \rangle=\left\{\left[\begin{matrix} 0&n \\ m&0 \end{matrix}\right], \left[\begin{matrix} m&0 \\ 0&n \end{matrix}\right]
\, :\, m,n \in \Z, mn\neq 0\right\}
$$
Every involution $I_r$  can be expressed in terms of the following list of involutions:
$$
V(z)=-z, \quad S(z)=I_{-1}(z)=-\frac{1}{z},\quad  I_{2}(z)=\frac{2}{z}, \quad I_{3}(z)=\frac{3}{z},\dots, I_p, \dots\mbox{ ($p$ prime).}
$$
For example, $I_{pq}=I_pI_1I_q=I_pUI_q$ where $U(z)=I_{1}=SV$.
The involution $V=SU=I_{-1}I_1$ commutes with all involutions $I_n$ for $n\in \Z$.
Also, $S=I_{-1}$ and $U=I_1$ commute, otherwise $I_p$ and $I_q$ do not commute.
\begin{proposition}
The set $\{K, S, V\}\cup \{I_p\, :\, p: \mbox{prime}\}$  generates $\pglq$.
\end{proposition}
{\it Proof}. Denote the translation $T_{r/s}(z):=z+r/s$. Then 
$T_{r/s}=H_{r/s}KVH_{s/r}$, and unless $p=0$, we have
$$
\frac{pz+q}{rz+s}=
\frac{p}{r}\left\{1-\frac{\frac{ps-qr}{pr}}{z+\frac{s}{r}}\right\}
=H_{p/r}KH_{(ps-qr)/pr}T_{s/r}(z).
$$
If $p=0$, then 
$$
\frac{pz+q}{rz+s}=\frac{q}{rz+s}=
\frac{q}{r}\frac{1}{z+\frac{s}{r}}
=H_{q/r}UT_{s/r}(z).
$$
Now   $U=KI_p(KV)^pVI_pK$ for any $p$, as one may easily check. 
Finally, the result follows since we can express the homotheties in terms of involutions. $\Box$
\begin{corollary}
The set $\{T, U, V\}\cup \{H_p\, :\, p: \mbox{ is prime }\}$  generates $\pglq$.
\end{corollary}
When $p=2$  the relation $U=KI_p(KV)^pVI_pK$ implies 
$U=KI_2KVKI_2K$, i.e. either one of the generators $V$ and $U$ can be eliminated from a generating set.
The case $p=2$ also implies that $U$ and $V$ are conjugate elements in  $\pglq$. 
They are not conjugates inside  $\pgl$.
\section{Elements of finite order}\label{finite}
$\PGL_2(\R)$ contains elements of any order. 
This is not true for $\pglq$.
Finite order matrices are elliptic so have $\mathrm{tr}^2(M)/\det(M)<4$.
It is known \cite{dresden} that the order of $M$ can be 2,3,4 or 6 and $M$ is $\pglq$-conjugate to
\begin{align}
\label{order3} M_3:=&\left[\begin{matrix} 0&-1 \\ 1&1 \end{matrix}\right]
\,\Bigl(\mbox{with fixed points } \frac{-1\pm i\sqrt{3}}{2}\Bigr)\, \mbox{ if the order is 3}\\
\label{order4} M_4:=&\left[\begin{matrix} 1&-1 \\ 1&1 \end{matrix}\right]
\,\Bigl(\mbox{with fixed points } \pm i\Bigr)\, \mbox{ if the order is 4}\\ 
\label{order5} M_6:=&\left[\begin{matrix} 2&-1 \\ 1&1 \end{matrix}\right]
\,\Bigl(\mbox{with fixed points } \frac{1\pm i\sqrt{3}}{2}\Bigr)\, \mbox{ if the order is 6}.
\end{align}
\paragraph{Elements of order two.}
$$
M=\left[\begin{matrix} p&q \\ r&s \end{matrix}\right] \implies M^2=
\left[\begin{matrix} p^2+qr&pq+qs \\ rp+rs&rq+s^2 \end{matrix}\right]=
\left[\begin{matrix} 1&0 \\ 0&1 \end{matrix}\right]
\implies 
$$
$$p^2+qr=s^2+rq \mbox{ and } q(p+s)=r(p+s)=0.
$$
Hence, $p^2=s^2$ and there are two possibilities: If $p=s\neq 0$ yields the identity.
Otherwise $p=-s$, and $q$, $r$ are free parameters, yielding the matrices of type
$$
M=\left[\begin{matrix} p&q \\ r&-p \end{matrix}\right], \mbox{ with } \det(M)=-p^2-qr\neq 0
$$
These account for the set of all non-singular traceless matrices. 
\paragraph{Elements of order three.} 
Routine calculations shows that if $M^3=I$ then
$$ 
M=\left[\begin{matrix} p/(p+s)&q/(p+s) \\ -(p^2+sp+s^2)/q(p+s)&s/(p+s) \end{matrix}\right]\in \pslq
\mbox{ if } p,q,s\in \Q
$$
We observe that these matrices are of trace one. 
\paragraph{Elements of order four.} Routine calculations shows that if $M^4=I$ then
$$
M=\left[\begin{matrix} p&q \\ -(p^2+s^2)/2q&s \end{matrix}\right] \quad q(p+s)\neq 0 \mbox{ with } \det =\frac{(p+s)^2}{2}
$$
These are not in $\pslq$. 
Their normalized trace is $\sqrt{2}$.
\paragraph{Elements of order six.} Routine calculations shows that if $M^6=I$ then
$$
M=\left[\begin{matrix} p&q \\ -(p^2-ps+s^2)/3q&s \end{matrix}\right] \quad q(p+s)\neq 0 \mbox{ with } \det =\frac{(p+s)^2}{3}
$$
These are not in $\pslq$. 
Their normalized trace is $\sqrt{3}$.

\medskip
\noindent
{\bf Remark:} 
One may ask the question, over which fields $K$, the group ${\mathrm{PGL}}(2,K)$ admits elements of  finite orders other then $2,3,4$ and $6$? 
To handle this general case we introduce a new coordinate system as follows: 
$$
\left[\begin{matrix} p&q \\r&s \end{matrix}\right]
=\left[\begin{matrix} x+y&z-t \\ z+t&x-y \end{matrix}\right], \quad {\left\{
\begin{array}{ll}
t:=\sqrt {{y}^{2}+{z}^{2}-{\delta}^{2}}\\ 
x^2/\delta^2:=\xi \end{array} \implies
\right.}
$$
$$
\left[\begin{matrix} x+y&z-t \\ z+t&x-y \end{matrix}\right]^n=
I
\iff
\left[
\begin{matrix} 
\delta\,\sqrt {\xi}+y&z-\sqrt {{y}^{2}+{z}^{2}-{\delta}^{2}}\\
z+\sqrt {{y}^{2}+{z}^{2}-{\delta}^{2}}&\delta\,\sqrt {\xi}-y
\end{matrix}
\right]^n=
I
$$
For $n$ odd, this reduces to just one equation of degree $(n-1)/2$ in $\xi$ which reads
$$
\begin{array}{ll}
n=3\implies& 3\,\xi+1 \\ 
n=5\implies& 5\,{\xi}^{2}+10\,\xi+1\\  
n=7\implies& 7\,{\xi}^{3}+35\,{\xi}^{2}+21\,\xi+1, \, etc.\\ 
\end{array}
$$
It can be shown that elements of higher orders appear in cyclotomic fields. The group 
$\mathrm{PGL}(2,K)$ contains elements of any finite order, if $K$ is the maximal abelian extension of $\Q$. 
We believe that for $K$ a number field, it should be possible to give a presentation of this group, in a way similar to we do here.

\section{Relators}
Here is our conjectural presentation for  $\pglq$:
\begin{center}
\fbox{\begin{minipage}{15cm}

{\bf Generators:} $\{T, U, V\!=\!H_{-1}, H_2,\dots H_p, \dots \, |\, p: \mbox{prime} \}$. \\
{\bf Relators:}
$$
\begin{array}{ll}
\mathbf{(I)}\quad U^2=V^2=(UV)^2=1 & \mathbf {(II)}\quad [H_p, H_q]=1\forall p, q\\
\mathbf{(III)} \quad(UH_p)^2=1\forall p &  \mathbf{(IV)}\quad H_p\inv T^pH_p=T\forall p \\
\mathbf{(V)} \quad V=T\inv UTUT\inv U&\mathbf{(VI)} \quad(TUT\inv U)^3=1\, (redundant)\\
\mathbf{(VII)} \quad(H_2UT\inv UT)^4=1&\mathbf {(VIII)}\quad (H_3UTUTUT^{-2})^6=1\\
\end{array}
$$
{\bf Dictionary:} \\
\noindent
\begin{tabular}{l}
$T(z):=z+1$: Translation\\
$V(z):=H_{-1}(z)=-z$: Reflection\\
$H_p(z)=pz$: Homothety\\
$TUT\inv U(z)=R^2(z)=1/(1-z)$: 3-rotation\\
$H_2UT\inv UT(z)=-2(z+1)/z$: 4-rotation\\
$H_3UTUTUT^{-2}(z)=(3z-3)/(2z-3)$: 6-rotation\\
\end{tabular}
\end{minipage}}
\end{center}

It is straightforward to check the validity of these relators. The difficulty lies in proving that there are no other relators independent from the above ones. This list have been found by a non-systematic search.
The last three equations originates from the finite order elements. The redundancy of the relator (VI) will become evident in the alternative presentation we give below.

\subsection{Another presentation of $\pglq$}
Here we rewrite the above presentation in terms of involutions, by expressing the homotheties in terms of $I_p$'s. The aim is to provide a framework for future research.

\begin{center}
\fbox{\begin{minipage}{15cm}
\medskip\noindent
{\bf Generators: } 
$\{K, V, U\!=\!I_{1}, I_2,\dots I_p, \dots \, |\, p: \mbox{prime} \}$.

\noindent
{\bf Relators:}
$$
\begin{array}{ll}
\mathbf{(I)}\quad U^2=V^2=(UV)^2=1 & \mathbf {(II)}\quad (I_pUI_q)^2=1\forall p, q\\
\mathbf{(III)} \quad I_p^2=1\forall p &  \mathbf{(IV)}\quad KI_p (KV)^pI_pVK=U \forall p \\
\mathbf{(V)} \quad (KU)^3=1&\mathbf{(VI)} \quad(KU)^6=1 \mbox{ (redundant)}\\
\mathbf{(VII)} \quad(I_2VKUV)^4=1&\mathbf {(VIII)}\quad (I_3UKVKUVK)^6=1\\
\end{array}
$$
\noindent
{\bf Dictionary:} \\
\begin{tabular}{l} 
$I_{-1}=S=UV=-1/z$\\
$L=KU =1-1/z$\\
$I_{1}I_{-1}=H_{-1}=V=-z$\\
$T=LS=KV=KI_{1}I_{-1}=z+1$\\
$I_{1}=U=SV=1/z$\\
$K=1-z$
\end{tabular}
\end{minipage}}
\end{center}

\bigskip\noindent
{\bf Acknowledgements.}
This research is funded by a Galatasaray University research grant and  T\"UB\.ITAK grant 115F412.

% ----------------------------------------------------------------

\end{document}